\documentclass[letterpaper,11pt]{article}
\usepackage{amsmath, amsthm, amssymb}    
\usepackage{bridges}			
\usepackage{graphicx}			
\usepackage[colorlinks=true, urlcolor=blue, citecolor=black, linkcolor=black]{hyperref}  
\usepackage{subfloat}
\usepackage{subcaption}			

\usepackage[utf8]{inputenc}
\usepackage{amsfonts}
\usepackage{mathtools}
\usepackage{wrapfig}
\usepackage{tikz}
\usepackage{listings}
\usepackage{comment}

\title{Oriented and Non-oriented Cubical Surfaces in The Penteract}

\author{Manuel Est\'evez\textsuperscript{1}, \'Erika Rold\'an\textsuperscript{1,2}, and Henry Segerman\textsuperscript{3}
\vspace{10pt}\\
\textsuperscript{1} ScaDS.AI, Leipzig University; estevez@mis.mpg.de\\
\textsuperscript{2} Max Planck Institute for Mathematics in the Sciences; roldan@mis.mpg.de\\
\textsuperscript{3} Oklahoma State University; henry@segerman.org
} 

\date{}					

\begin{document}

\maketitle

\thispagestyle{empty}

\begin{abstract}
 Which surfaces can find the two-dimensional faces of the five-dimensional cube (the penteract)? How can we visualize them? In recent work, Aveni, Govc, and Rold\'an show that there exist 2690 connected closed cubical surfaces up to isomorphism in the 5-cube. They give a classification in terms of their genus $g$ for closed orientable cubical surfaces, and their demigenus $k$ for a closed non-orientable cubical surface. In this paper we explain the main idea behind the exhaustive search and we visualize the projection to $\mathbb{R}^3$ of a torus, a genus two torus, the projective plane, the Klein bottle. We use reinforcement learning techniques to obtain configurations optimized for 3D-printing.
\end{abstract}
\section{Introduction}
In Bridges 2023 \cite{EstevezRoldanSegerman}, Est\'evez, Rold\'an, and Segerman presented various visualizations of 3D-printed  representations of \textit{closed connected orientable cubical surfaces} embedded in $\mathbb{R}^3$ on the two dimensional faces of the tesseract. These representatives are homeomorphic to either a sphere, a torus or two disconnected spheres; in fact, the maximum genus achievable for the tesseract is one. For a five-dimensional cube, often called a \textit{penteract}, doing an exhaustive computational search, Aveni, Govc, and Rold\'an \cite{AveniGovcRoldan} found all orientable and non-orientable connected closed surfaces. Within these surfaces, by a result from Schulz \cite{Schulz}, the maximum possible genus is five and the maximum possible demigenus is eight. They have also classified all closed surfaces in 2690 different isomorphism types. 

Here, we have selected some of these surfaces to find good embeddings for visualizing them in $\mathbb{R}^3$ and for 3D-printing them: the torus, the projective plane, and the Klein Bottle. We work with non-orientable cubical surfaces, therefore we must deal with self intersections of its faces; in particular we want to minimize them. In order to have the best possible embedding for 3D printing and visualization of these surfaces, we implement a Reinforcement Learning (RL) algorithm that explores suitable three-dimensional embeddings which minimize self intersections of its faces and edge crossings for a fixed edge width.

The rest of the paper is organized as follows. In Section \ref{CubicalSurfaces}, we introduce the notation and basic definitions of cubical surfaces in the penteract. In Section \ref{Visualization} we present some results and 3D-models of the configurations that the algorithm found for our selected surfaces, starting from a particular initial configuration. In Section \ref{RL}, we describe in detail the implementation of the RL algorithm.

\section{Cubical Surfaces in the Penteract}\label{CubicalSurfaces}
We denote the five-dimensional unit cube by $Q^5=[0,1]^5$, and its set of vertices by $Q_0^5$. Each vertex of $Q^5$ can be represented by an element of the set of all five-tuples with binary entries $\{0,1\}^5$. We denote by $Q_1^5$ the one-dimensional skeleton of $Q^5$, that is, the set of its vertices $v$ and edges $e$. We observe that $Q_1^5$ is the graph with vertex set $Q_0^5$ and an edge between two vertices if and only if they differ in exactly one coordinate. The cardinality of the set of faces $f$ containing an edge $e$ (resp. a vertex $v$) is denoted by $F_e$ (resp. $F_v$) and similarly the cardinality of the set of edges $e$ containing a vertex $v$ is denoted by $E_v$. Similarly, $Q_2^5$ denotes the two-dimensional skeleton of $Q^5$, which consists of the set of vertices $Q^5_0$, the one-dimensional skeleton $Q^5_1$, and all its two-dimensional faces $f$. We can continue this construction up to the penteract $Q^5_5$ itself, and name the elements of all the preceding sets the \textit{cells} of $Q^5$. Geometrically, every cell of $Q^5$ is a product of vertices and intervals, and therefore can be encoded combinatorially as an element of $\{0,1,*\}^5$. Here a * in an entry implies that in the product, the whole interval $I$ is considered in that direction. Thus, every sub-complex of $Q^5$ can be represented as a subset of $\{0,1,*\}^5$. We call this, the \textit{star notation}. We refer to a subset of $Q^n_2$ as a \emph{two-dimensional cubical complex}, which we will denote by $\mathcal{C}$. We denote its set of vertices, edges and faces by $\mathcal{C}_0$, $\mathcal{C}_1$, and $\mathcal{C}_2$ respectively. The \textit{vertex figure} $\mathcal{F}_v$ of a vertex $v$ is the graph whose nodes are the edges in $\mathcal{C}_1$ having $v$ as an endpoint and where two nodes $e,e' \in C_1$ are joined by an edge if and only if there is a face $f \in \mathcal{C}_2$ with $e,e'$ as two of its edges. A \textit{closed cubical surface} is a two-dimensional cubical complex $\mathcal{C}$ in which every point has an open neighborhood homeomorphic to an open disk. This condition is equivalent to asking $\mathcal{C}$ to fulfill the following conditions:
\begin{enumerate}
\item Every edge is shared by exactly two faces, i.e. $\forall e \in \mathcal{C}_1, F_e=2$ ;
\item The vertex figure $\mathcal{F}_v$ of any vertex $v \in \mathcal{C}_0$ is a cyclic graph.
\end{enumerate}
\section{Embeddings of Cubical Surfaces}\label{Visualization}
For the following surfaces we are using the results found by Aveni, Govc and Rold\'an \cite{AveniGovcRoldan}, therefore we know that we are using the minimum number of faces needed for realizing these orientable and non-orientable surfaces.
\subsection{Orientable Surfaces}
\begin{figure}[htbp]
    \centering
    \begin{subfigure}[b]{0.24\textwidth}
        \includegraphics[width=4cm]{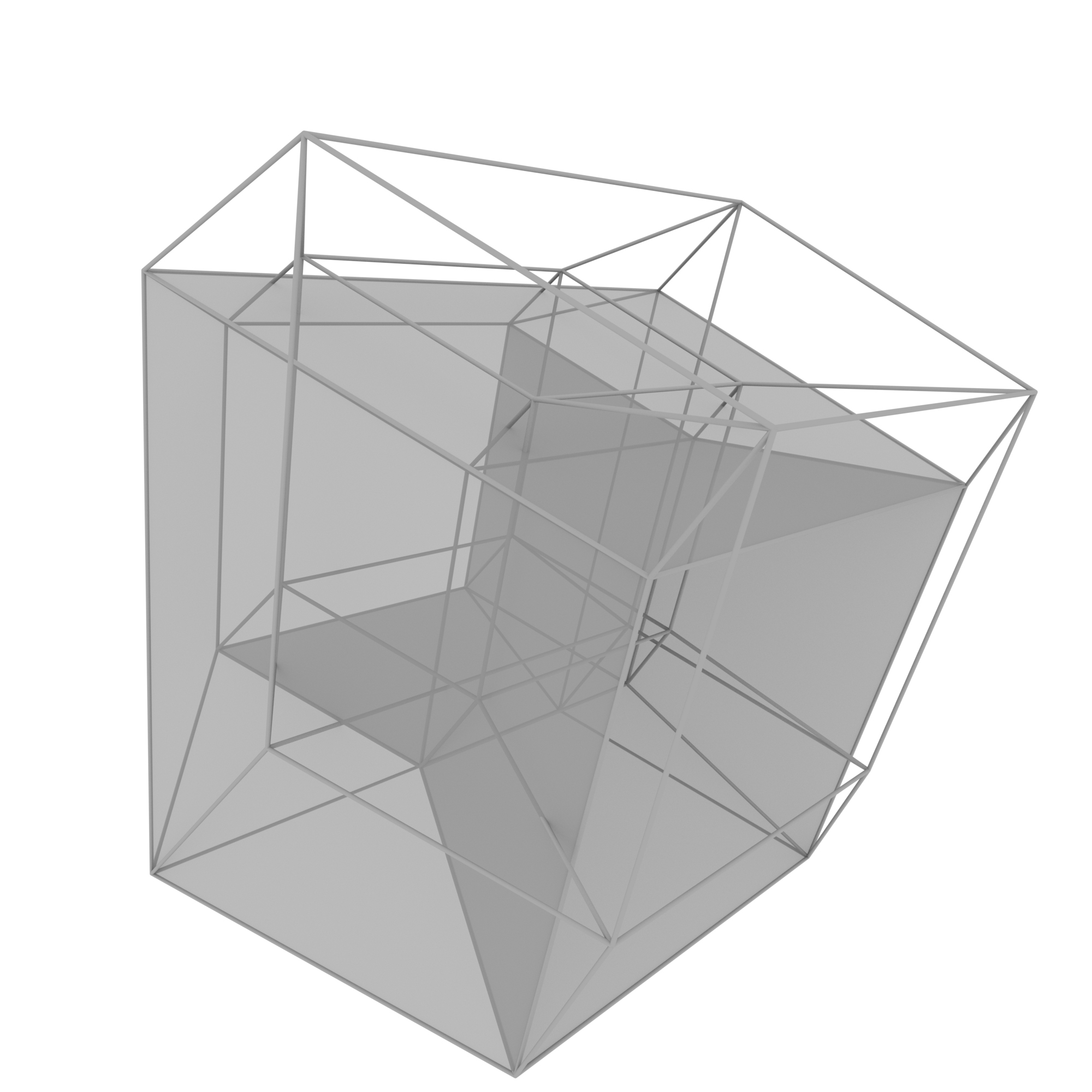}
        \caption{Genus 1 torus}
        \label{rfidtest_xaxis}
    \end{subfigure}
    \begin{subfigure}[b]{0.24\textwidth}
        \includegraphics[width=4cm]{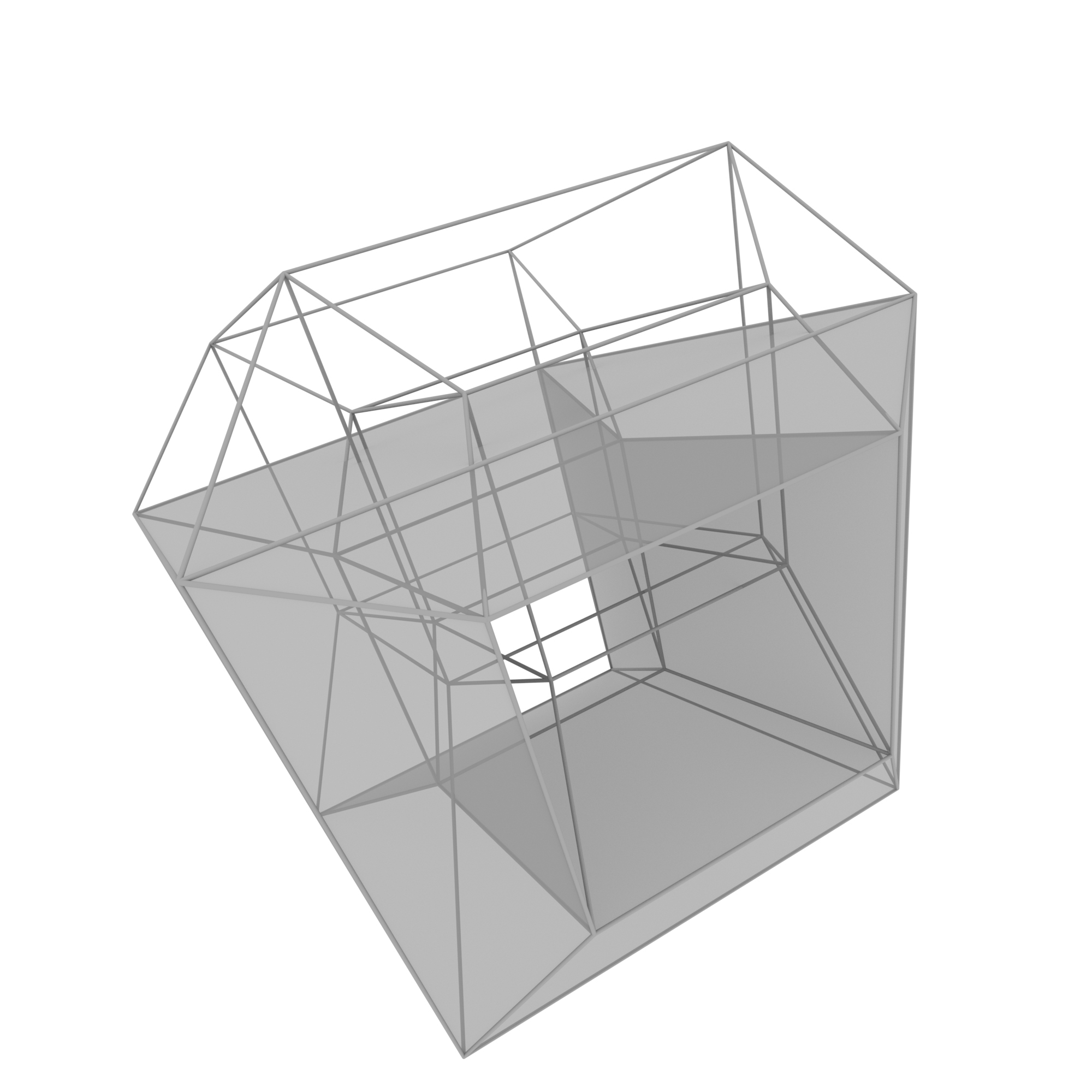}
        \caption{Edge overlaps 0 \\
        Face intersections 0}
        \label{rfidtest_yaxis}
    \end{subfigure}
    \begin{subfigure}[b]{0.24\textwidth}
        \includegraphics[width=4cm]{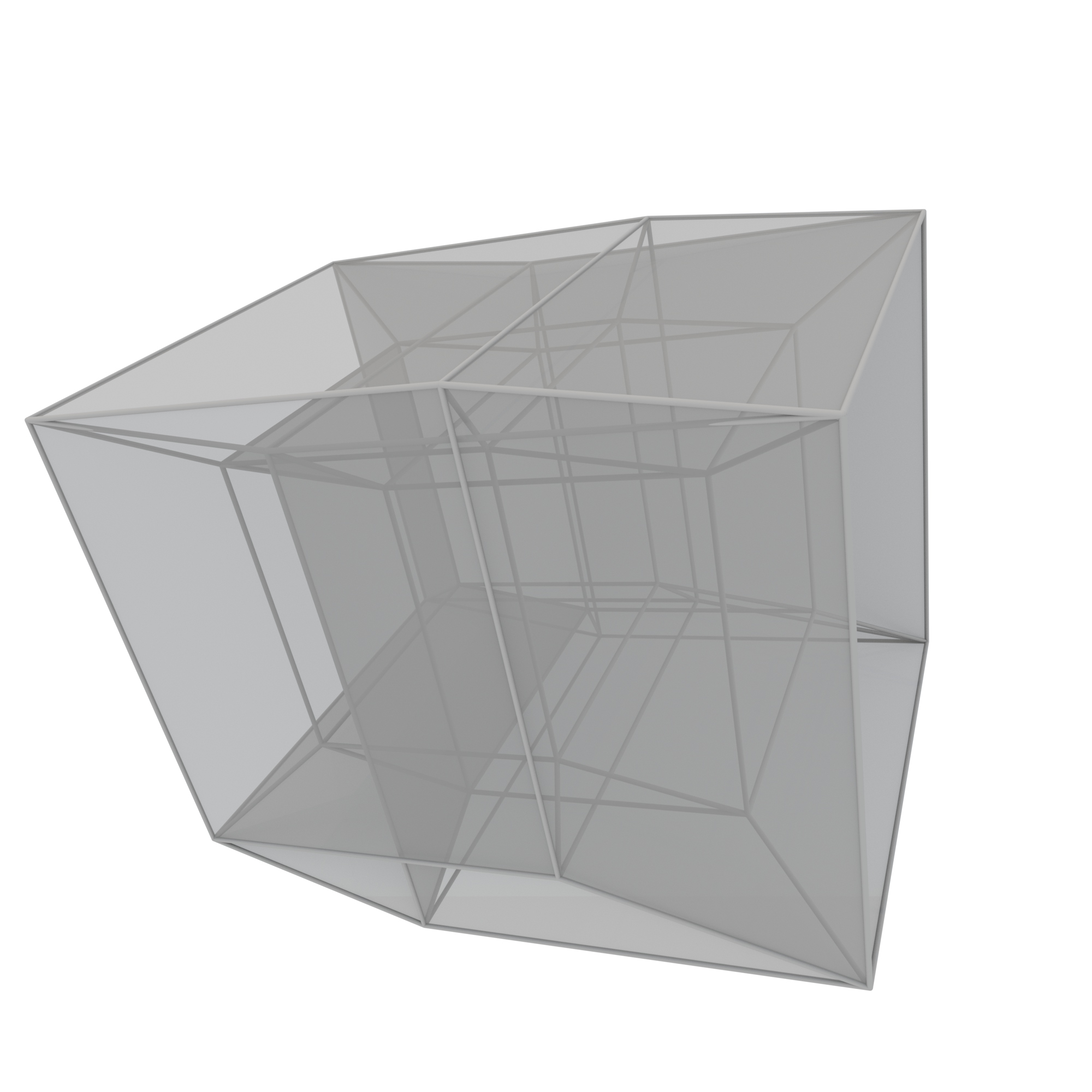}
        \caption{Genus 2 torus}
        \label{rfidtest_yaxis}
    \end{subfigure}
        \begin{subfigure}[b]{0.24\textwidth}
        \includegraphics[width=4cm]{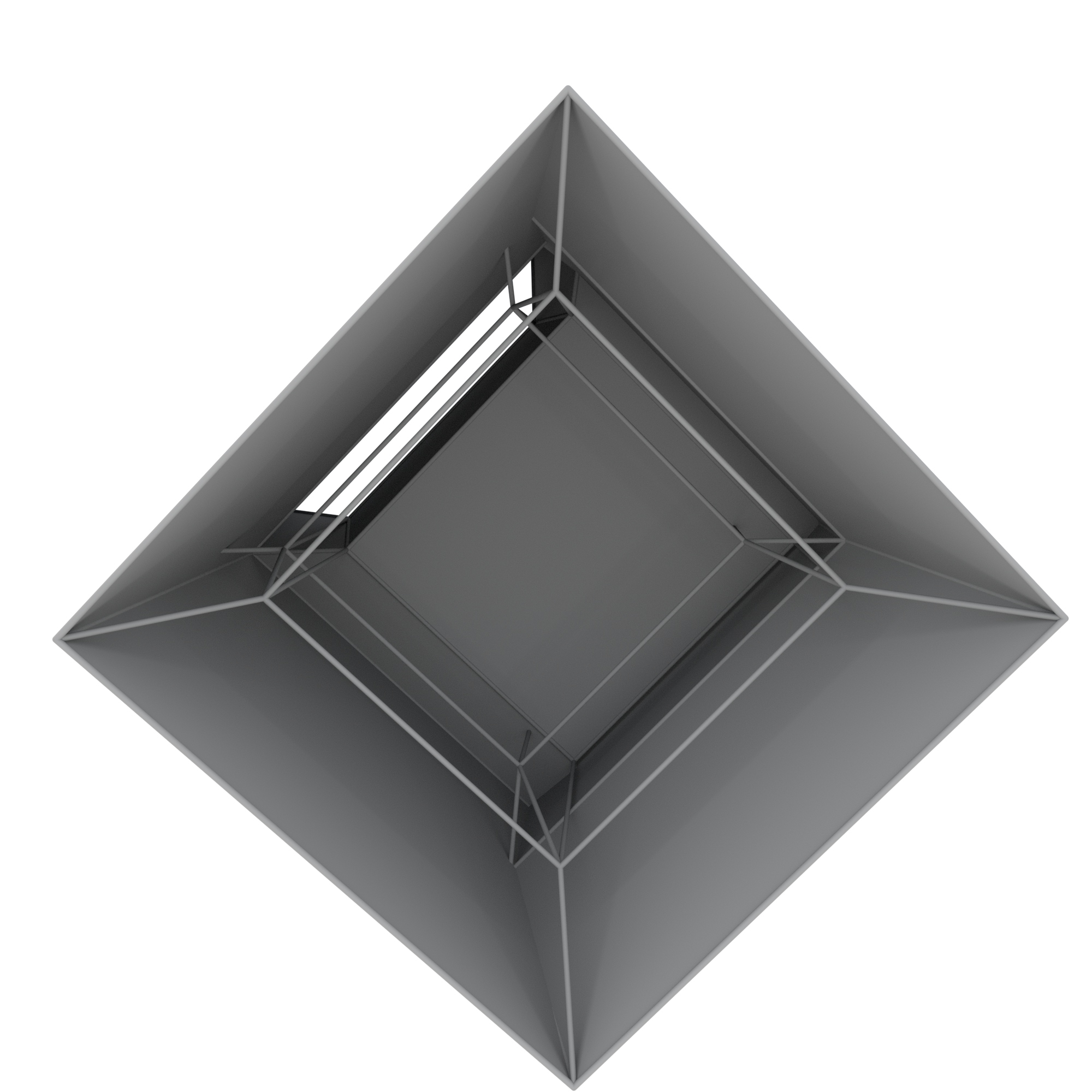}
        \caption{Edge overlaps 0 \\
        Face intersections 11}
        \label{rfidtest_yaxis}
    \end{subfigure}
    \caption[RFID tag read-range testing]{Orientable cubical surfaces} \label{orientablefigures}
\end{figure}

In Figure \ref{orientablefigures} (a) and (b) (resp. (c) and (d)) we present an embedding of a torus (resp. genus two torus) whose initial embedding had 19 (resp. 19) edge overlaps and 6 (resp. 33) face intersections. Our RL algorithm found embeddings with zero face overlaps for both and zero (resp. 11) face intersections.
\subsection{Non-orientable Surfaces}
\begin{figure}[htbp]
    \centering
    \begin{subfigure}[b]{0.32\textwidth}
        \includegraphics[width=5cm]{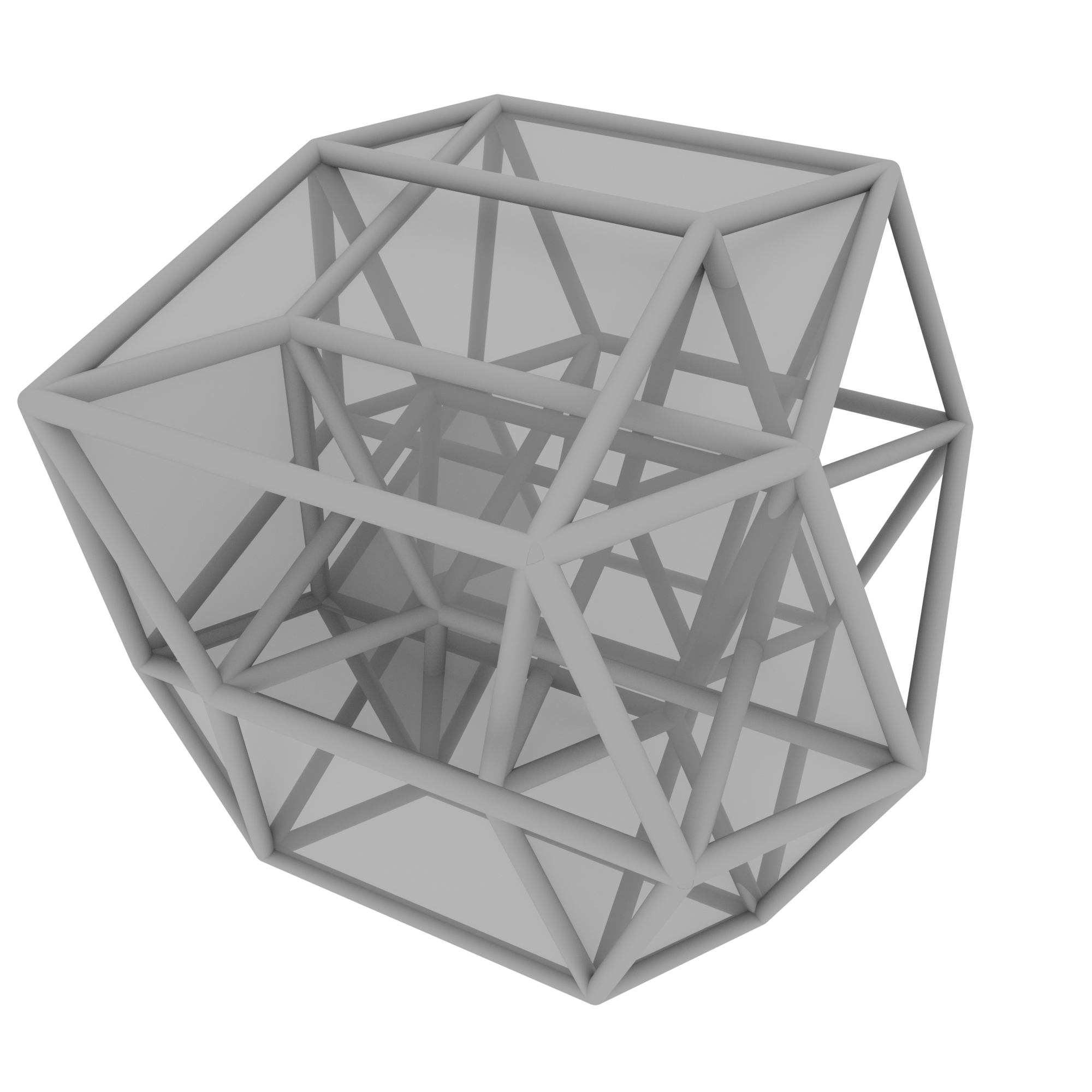}
        \caption{Edge overlaps 0 \\
        Face intersections 3}
        \label{rfidtest_yaxis}
    \end{subfigure}
    \begin{subfigure}[b]{0.32\textwidth}
        \includegraphics[width=5cm]{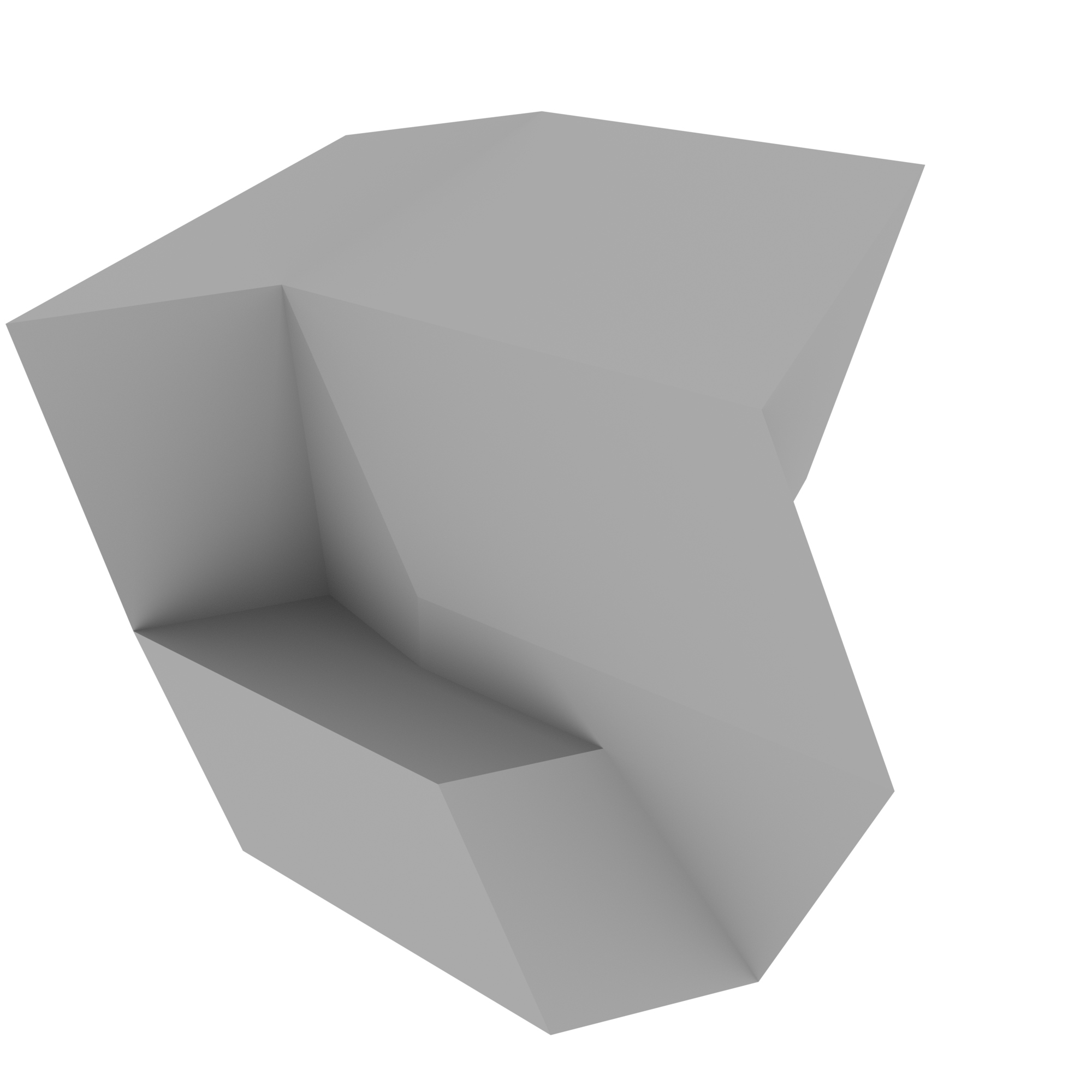}
        \caption{Projective plane's faces (front)}
        \label{rfidtest_yaxis}
    \end{subfigure}
    \begin{subfigure}[b]{0.32\textwidth}
        \includegraphics[width=5cm]{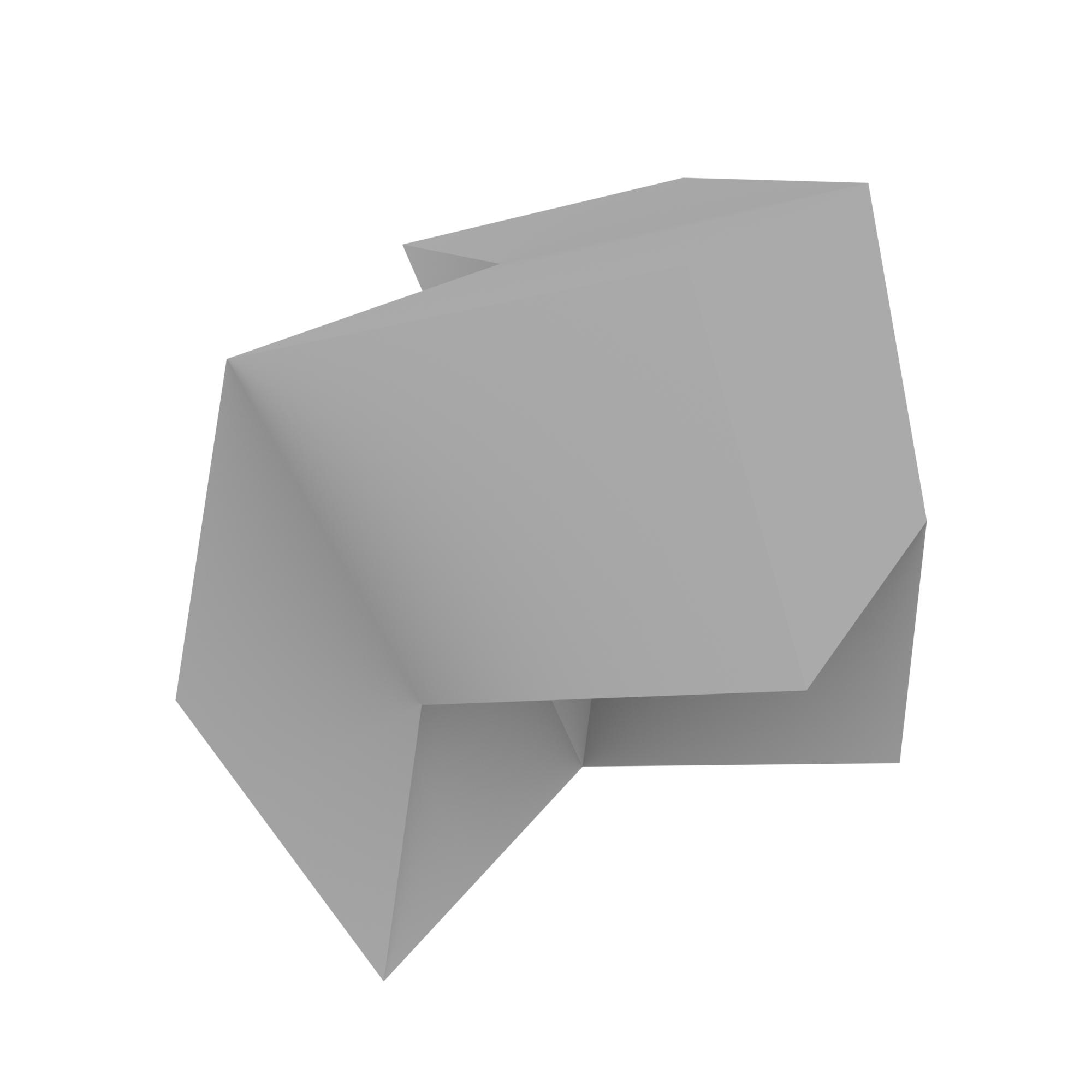}
        \caption{Projective plane's faces (back)}
        \label{rfidtest_yaxis}
    \end{subfigure}
    \caption[RFID tag read-range testing]{A cubical surface homeomorphic to the projective plane.}\label{Projplane}
\end{figure}

\begin{figure}[h]
    \centering
    \begin{subfigure}[b]{0.32\textwidth}
        \includegraphics[width=5cm]{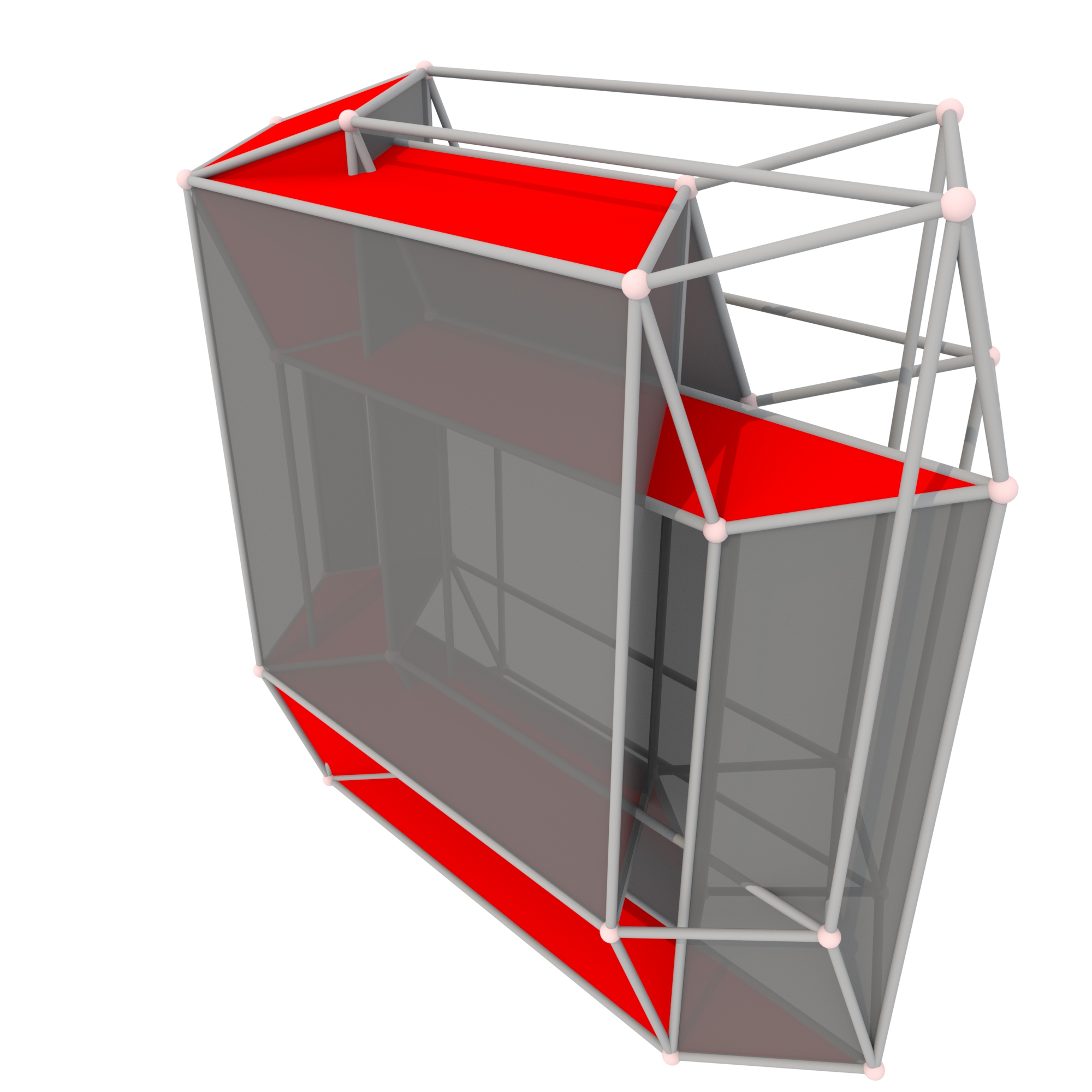}
        \caption{Edge overlaps: 0 \\
        Face intersections $>3$.}\label{Klein_multiple_intersections}
        \label{rfidtest_xaxis}
    \end{subfigure}
    \begin{subfigure}[b]{0.32\textwidth}
        \includegraphics[width=5cm]{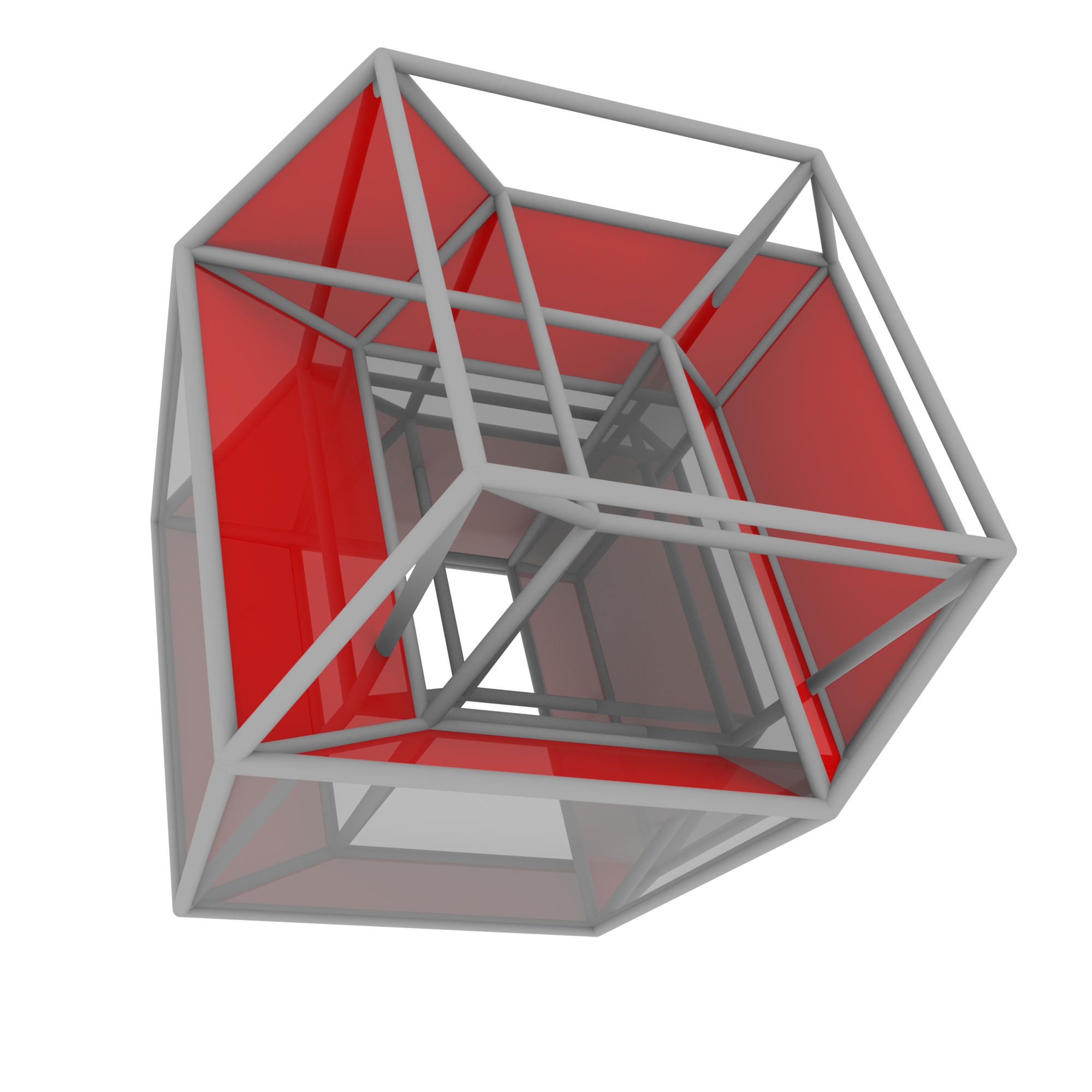}
        \caption{Klein Bottle with one-skeleton $Q^{5}_1$}
        \label{rfidtest_xaxis}
    \end{subfigure}
    \begin{subfigure}[b]{0.32\textwidth}
        \includegraphics[width=5cm]{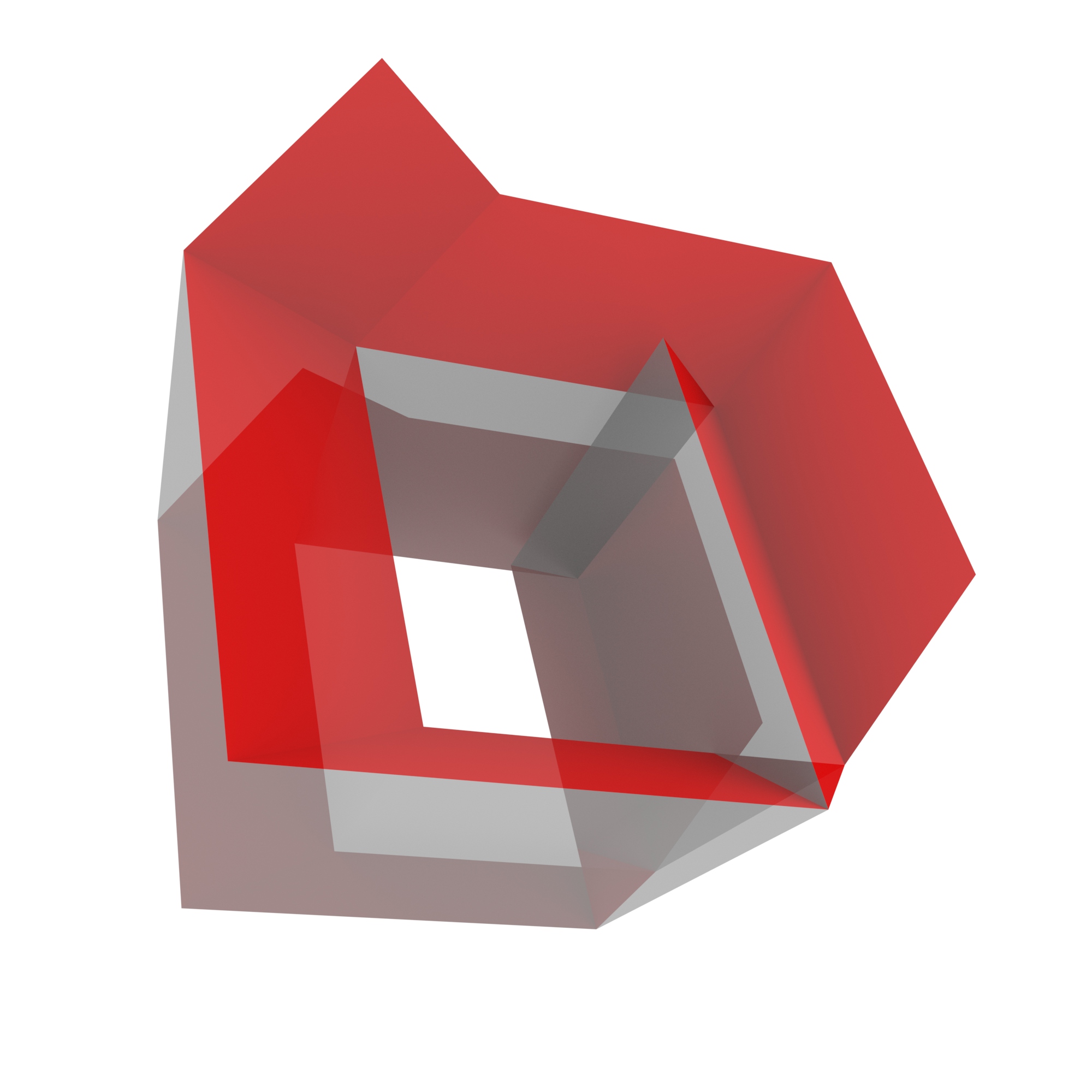}
        \caption{Klein Bottle with two M\"obius strips in red}
        \label{rfidtest_yaxis}
    \end{subfigure}
    \caption[RFID tag read-range testing]{A cubical surface homeomorphic to the Klein Bottle}\label{KB}
\end{figure}

In Figure \ref{Projplane} (resp. Figure \ref{KB}), we present an embedding of a cubical surface homeomorphic to the projective plane (resp. Klein Bottle), whose initial embedding had 19 (resp. 19) edge overlaps and 13 (resp. 9) face intersections. Our RL algorithm lowered the intersections to only three face intersections and zero edge overlaps for both of them. In Figure \ref{Klein_multiple_intersections}, we give an embedding of the Klein Bottle with a lot of intersections, emphasizing two M\"obius strips that are subcomplexes of this surface. We challenge the reader to count the number of intersections using our 3D model \url{https://skfb.ly/oRB7H}. 

Our projective plane and Klein bottle each have two cross-cap singularities. We are currently  working on modifying our code to search for immersions without these singularities.

\section{Using Reinforcement Learning to Minimize Face Intersections and Edge Overlaps} \label{RL}
The embedding of $\mathcal{C}$ in $\mathbb{R}^3$ at time step $t$ is parameterized by a vector $s_t=(d_5, d_4, \phi_1, \cdots , \phi_{10}) \in \mathbb{R}^{12}$ which we call a \textit{state}. The value $d_5$ (resp. $d_4$) is the distance from the camera point to the origin on perspective projection $pr_5(d_5): \mathbb{R}^5 \rightarrow \mathbb{R}^4$ (resp. $pr_4(d_4): \mathbb{R}^4 \rightarrow \mathbb{R}^3$) in five (resp. four)-dimensional space. When using perspective projection in five (resp. four)-dimensional space, the distance from the origin to the projection hyperplane must be specified, we fix these distance to be $c_5=1$ (resp. $c_4=10$). One the other hand, the values $\phi_i \in [0, 2\pi), (1 \leq i \leq 10)$ are the ten possible angles (one for each pair of axes) on which the penteract $Q^5$ and $\mathcal{C}$ can be rotated around the origin in $\mathbb{R}^5$. We call the set $\mathcal{S}$ of states, the \textit{state space}, and the initial state by $s_0\in \mathcal{S}$.  The RL algorithm to which we refer simply as the \textit{agent} can perform one of the possible actions in the \textit{action set} $\mathcal{A}:=\{ \delta e_0, -\delta e_0, \delta e_1, -\delta e_1, \epsilon e_2, -\epsilon e_2, \cdots , \epsilon e_{11}, -\epsilon e_{11} \}$, where the $e_i \in \mathbb{R}^{12}$ are vectors whose $i$ coordinate is 1 and $0$ elsewhere, and $\delta, \epsilon \in \mathbb{R}$ are positive real numbers defining the step length which we set after some tests to $\delta = .5$ and $\epsilon = \pi/180$.  An action $a_i \in \mathcal{A}$ is also a vector in $\mathbb{R}^{12}$, and given a state $s_t$ we can think on the state $s_{t+1}$ resulting by taking some action $a_i$ as the sum $s_{t+1}=s_{t}+a_i$. The parameter $\delta$ affects the distances $d_5$ and $d_4$ while the parameter $\epsilon$ affects the rotation angles $\phi_i, (1 \leq i \leq 10)$. Therefore, an action can apply either a small five-dimensional rotation or move the five or four- dimensional cameras to obtain a better embedding.  A \textit{reward function} $R: \mathcal{S} \times \mathcal{A} \rightarrow \mathbb{R}$ gives feedback to the agent depending on whether the action taken favors a certain task. We explain the construction of our reward functions in Sections \ref{Faces} and \ref{Edges} for different tasks. A \textit{transition distribution}  $P:\mathcal{S} \times \mathcal{A} \rightarrow (\mathcal{S} \rightarrow [0,1])$ determines the probability of arriving to a certain state $s_t$ given a tuple $(s_{t-1}, a_i) \in \mathcal{S} \times \mathcal{A}$. The tuple $M = ( \mathcal{S}, \mathcal{A}, R, P, \gamma, \rho_0)$, where $\gamma \in [0,1)$ is a discount factor and $\rho_0: \mathcal{S} \rightarrow [0,1]$ an initial state probability distribution is called a Markov Decision Process. Solving $M$ means finding a probability distribution $\Pi$ over $\mathcal{A}$ at each state that yields the supremum of the expected discounted sum of rewards $\underset{\Pi}{\mathrm{sup}}\, \mathbb{E}_{s_0\sim \rho_0} \left[\mathbb{E}_\Pi\left[\sum_{t=0}^{\infty} \gamma^t R_{t+1}|S_0 = s_0\right]\right]$; $P$ tends to $\Pi$ as $t \rightarrow \infty$.

\subsection{Minimizing the number of intersecting faces} \label{Faces}
For a pair of distinct faces $(f_i,f_j) \in \binom{\mathcal{C}_2}{2}$ after performing some five-dimensional rotation, consider their perspective projections $pr(d_5, d_4,f_i)$ and $pr(d_5, d_4,f_j)$ in $\mathbb{R}^3$. If these projections intersect, they do it transversely (in a point or a line) or they overlap (we don't consider faces sharing an edge as intersecting). At the state $s_t \in \mathcal{S}$ the \textit{number of face intersections} running through all pairs $(i,j)$ is denoted by $\Sigma_{\mathcal{C}}(s_t)$.  For most of our cubical surfaces we don't know the minimum number of face intersections, but we can propose a minimum ${\Sigma_{\mathcal{C}}}_{prop}$ and expect the algorithm to find a state $s_t$ such that $\Sigma_{\mathcal{C}}(s_t) \leq {\Sigma_{\mathcal{C}}}_{prop}$.  The reward function
\begin{equation*}
R_1(s_t, a_t) := \begin{cases}
10*(1-\Sigma_{\mathcal{C}}(s_t+a_t)+{\Sigma_{\mathcal{C}}}_{prop}) &\text{if $\Sigma_{\mathcal{C}}(s_t+a_t) \leq {\Sigma_{\mathcal{C}}}_{prop}$}\\
\frac{\Sigma_{\mathcal{C}}(s_t+a_t)-\Sigma_{\mathcal{C}}(s_t)}{\Sigma_{\mathcal{C}}(s_t)} &\text{if $\Sigma_{\mathcal{C}}(s_t+a_t) > {\Sigma_{\mathcal{C}}}_{prop}$},
\end{cases}
\end{equation*}
will reward the agent when the action $a_t$ reduces $\Sigma_{\mathcal{C}}(s_t+a_t)$ with respect to $\Sigma_{\mathcal{C}}(s_t)$.
\subsection{Minimizing edge overlaps for 3D-printing}\label{Edges}
For a pair of distinct edges $(e_i,e_j) \in \binom{\mathcal{C}_1}{2}$ after performing some five-dimensional rotation, consider their perspective projections $pr(d_5,d_4,e_i)$ and $pr(d_5,d_4,e_j)$ in $\mathbb{R}^3$. To 3D-print the three-dimensional projection of the cubical surface we assign a constant width $r$ to all of the projected edges. Each pair of projected edges can overlap with each other at a given state $s_t \in \mathcal{S}$ if the perpendicular line segment $L_{e_i,e_j}(s_t)$ connecting them has magnitude $|L_{e,e'}(s_t)|<2r$. The \textit{number of overlapping edges} at a state $s_t \in \mathcal{S}$ is denoted by $o_w(s_{t})$. The reward function
\begin{equation*}
R_2(s_t, a_t) := \begin{cases}
10*(1-\Sigma_{\mathcal{C}}(s_t+a_t)+{\Sigma_{\mathcal{C}}}_{prop}) &\text{if $o(s_{t})=0$ and $\Sigma_{\mathcal{C}}(s_t+a_t) \leq {\Sigma_{\mathcal{C}}}_{prop}$}\\
\frac{o(s_t+a_t)-o(s_{t})}{o(s_{t})} &\text{if $o(s_{t}) \neq 0$ and $\Sigma_{\mathcal{C}}(s_t+a_t) \leq {\Sigma_{\mathcal{C}}}_{prop}$}\\
0 &\text{if $\Sigma_{\mathcal{C}}(s_t+a_t) > {\Sigma_{\mathcal{C}}}_{prop}$},
\end{cases}
\end{equation*}
will reward the algorithm when the action $a_t$ reduces $o_w(s_t+a_t)$ with respect to $o_w(s_0)$ if the number ${\Sigma_{\mathcal{C}}}_{prop}$ is achieved or improved. We prevent the agent from lowering $o(s_{t})$ simply by decreasing the parameters $d_5, d_4$ because we want the width $r $ not to be too small with respect to the final size of the projection for 3D printing purposes. At a state $s_t \in \mathcal{S}$, we consider $L(s_t):= \sum_{(e_i,e_j) \in \binom{\mathcal{C}_1}{2}}|L_{e_i,e_j}(s_t)|$. The function $R_3(s_t, a_t) := \frac{L(s_t)-L(s_t+a_t)}{L(s_t)}$, will reward the agent for reducing $L(s_t)$. We take $R_4(s_t, a_t):=1$ if $L(s_t)< min(\{L(s_i)_{0 \leq i \leq t-1}\})$ and $R_4(s_t, a_t):=0$ otherwise. We take then $R:=R_1+R_2+R_3+R_4$.

\newpage
{\setlength{\baselineskip}{13pt} 
\raggedright				

} 
\end{document}